# Symmetric similarity 3D coordinate transformation based on dual quaternion algorithm

Sebahattin Bektas

Ondokuz Mayis University,  Geomatics Engineering,  Samsun, Turkey,

sbektas@omu.edu.tr

ORCID: 0000-0001-8210-0352

## ABSTRACT

Nowadays, we have seen that dual quaternion algorithms are used in 3D coordinate transformation problems due to their advantages. 3D coordinate transformation problem is one of the important problems in geodesy. This transformation problem is encountered in many application areas other than geodesy. Although there are many coordinate transformation methods (similarity, affine, projective, etc.), similarity transformation is used because of its simplicity. The asymmetric transformation is preferred to the symmetric coordinate transformation because of its ease of use. In terms of error theory, the symmetric transformation should be preferred. In this study, the topic of symmetric similarity 3D coordinate transformation based on the dual quaternion algorithm is discussed, and the bottlenecks encountered in solving the problem and the solution method are discussed. A new iterative algorithm based on the dual quaternion is presented. The solution is implemented in two different models: with constraint equations and without constraint equations. The advantages and disadvantages of the two models compared to each other are also evaluated. Not only the transformation parameters but also the errors of the transformation parameters are determined. The detailed derivation of the formulas for estimating the symmetric similarity of 3D transformation parameters is presented step by step. Since symmetric transformation is the general form of asymmetric transformation, we can also obtain asymmetric transformation results with a simple modification of the model we developed for symmetric transformation. The proposed algorithm is capable of performing both 2D and 3D symmetric and asymmetric similarity transformations. For the 2D transformation, it is sufficient to replace the z and Z coordinates in both systems with zero.

## 1. Introduction

Recently, we have seen many coordinate transformation applications based on dual quaternion algorithms (DQA). The 3D coordinate transformation problem is also encountered in GNSS and Lidar applications Li et al. (2013). As is known, unlike Euler angles, quaternions determine the position of an object in space in a univocal way. Due to the advantages of quaternions, their use is increasing day by day. Bektas (2024b) used quaternions in the orthogonal ellipsoid fitting problem. For more information on the advantages of 3D transformations using quaternions, see Zeng et al. (2018), Uygur et al. (2022), Bektas (2022), Bektas (2024a), and Zeng et al. (2024). There are fewer studies on symmetric transformations in the literature, some of which are Zeng and Yi (2011), Fang (2015), Felus and Burtch (2009). Mahboub (2016), Mercan et al. (2018), Wang et al. (2023) and Zhao et al. (2024) are on symmetric transformation with quaternion algorithms (QA) method. In the literature, there is only Zeng et al. (2024) on symmetric transformation with (DQA) method.

The date of precision calculations in DQA is new. The precision calculations in asymmetric DQA-based transformations were first made by Bektas (2024a). More recently, precision calculations for the DQA-based symmetric transformation were given by Zeng et al. (2024).

As it is known, while the translation parameters are determined classically, the scale factor and rotation angles are determined from quaternions, in the quaternion algorithm QA method. The QA method is also called the single quaternion method. In the DQA method, all transformation parameters are determined from quaternions. There are several studies on quaternion-based 3D coordinate transformations in the literature. Some of them are QA-based, while some are DQA-based. In general, however, asymmetric transformation methods have always been used because of their simplicity. As is well known, asymmetric transformation assumes that only the coordinates of the second system are erroneous and the coordinates of the first system are error-free, which is inconsistent with reality. This assumption of asymmetric transformation is not correct from the point of view of error theory. In the symmetric transformation model, the coordinates of both systems are assumed to be erroneous. The symmetric transformation model also takes into account that the coordinates of both systems may have different weights and even correlations.

There are some studies in the literature showing that DQA, QA and Euler angle methods give different results. However, both the DQA and QA methods are derived from Helmert's seven-parameter similarity transform. Therefore, DQA and QA and even Euler Angles 3D transformation methods should give the same results regardless of rounding errors, as reported Bektas (2022) and Bektas (2024a). The confusion here needs to be resolved. For example, Ionnadis and Pantazis (2020) used a model with nine parameters (eight quaternions + one scale). It is assumed that they use unit quaternions in that study form, but additional constraint equations would have to be found for unit quaternions to exist. However, no additional constraints are mentioned in the study. On the other hand, the results of that study show that DQA, QA and Euler angle methods produce parameters with different precision. Bektas (2024c) has addressed the inconsistencies in this article.

Zeng et al. (2024) performed DQA-based symmetric 3D coordinate transformation using the Total Least Squares (TLS) method in their study. The authors also claim that the DQA method gives more precise results than QA in its solution. They found the same transformation parameters as DQA and QA methods (Actual geodetic datum transformation case, Page 13, Table 6,7, and 8). However, they claim that the precision of the transformation parameters found by the DQA and QA methods in Table 7 are significantly different and that the DQA method produces parameters with higher precision. The authors have made a comparison between DQA and QA methods using scaled quaternions in Table 7. However, the precision values of *$r_4$ ($q_0$), $t_x$, $t_y$, $t_z$* in the DQA column of Table 7 are incorrect. Unfortunately, the authors concluded that the DQA method gives more precise results.

The rest of the article is organised as follows: Section 2 provides brief information about quaternions and dual quaternions, this section presents the mathematical model of a new 3D symmetric transformation based on dual quaternion with constraint equations and unconstrained model, also in this section sensitivity calculations, as well as the mathematical model according to the QA method. Section 3 contains two numerical experiments are designed to verify the correctness and effectiveness of the two solutions. Finally, the results and experiment are presented in the conclusion.

## 2. Material and Methods

In this section, mathematical models of various quaternion-based symmetric 3D coordinate transformations will be introduced. The 3D transformation models based on DQA with constraint equations and unconstrained(simplify) model, QA model and using scaled quaternions.

### Quaternion and dual quaternion

The history of quaternions dates back to 1853. The quaternions were invented by Sir William Hamilton in 1843. The word quaternion comes from the meaning of quaternary and is usually explained as follows.

*q = q₁ i + q₂ j + q₃ k + q₄*  (1)

where $q_1, q_2, q_3, q_4$ all of real numbers $i, j$ and $k$ imaginary units and they have the following features

$$i^2 = j^2 = k^2 = -1 \tag{2}$$

$$ij = -ji = k, \quad jk = kj = i, \quad ki = -ik = j, \quad ijk = -1 \tag{3}$$

the norm of quaternion

$$\|q\| = \sqrt{q_1^2 + q_2^2 + q_3^2 + q_4^2} \tag{4}$$

if $\|q\| = 1$, $q$ is named unit quaternion

The rotation matrix $R$ can be defined in terms of quaternions $q$ as follows

$$R = (q_4^2 - q^T q) I_{3x3} + 2(qq^T + q_4 C(q)) \tag{5}$$

$$q = [\ q_1\ q_2\ q_3\ ]^T \qquad r = [\ r_1\ r_2\ r_3\ ]^T \tag{6}$$

where $I_{3x3}$ is the unit matrix and $C(r)$ is a skew-symmetric matrix

$$C(r) = \begin{bmatrix} 0 & -r_3 & r_2 \\ r_3 & 0 & -r_1 \\ -r_2 & r_1 & 0 \end{bmatrix} \tag{7}$$

The rotation matrix $R$ consists of unit quaternions

$$R = \begin{bmatrix} q_4^2 + q_1^2 - q_2^2 - q_3^2 & 2(q_1 q_2 - q_4 q_3) & 2(q_1 q_3 + q_4 q_2) \\ 2(q_1 q_2 + q_4 q_3) & q_4^2 - q_1^2 + q_2^2 - q_3^2 & 2(q_2 q_3 - q_4 q_1) \\ 2(q_1 q_3 - q_4 q_2) & 2(q_2 q_3 + q_4 q_1) & q_4^2 - q_1^2 - q_2^2 + q_3^2 \end{bmatrix} \tag{8}$$

The dual quaternion which was invented by Clifford in 1873. Clifford (2007) first demonstrated dual quaternion for rotations and translations in a single model (Tucker,1968).

$$q = r + s.\mu \tag{9}$$

$\mu$ is a dual unit with the property $\mu^2 = 0$ and $\mu$ commutes with quaternion units (Zeng et al., 2018). Where $r$ and $s$ are quaternions. To perform 3D coordinate transformation with dual quaternions, the following matrix definitions are made.

$$W_{(r)}=\begin{bmatrix} r_4 I - C(r) & r \\ -r^T & r_4 \end{bmatrix} = \begin{bmatrix} r_4 & r_3 & -r_2 & r_1 \\ -r_3 & r_4 & r_1 & r_2 \\ r_2 & -r_1 & r_4 & r_3 \\ -r_1 & -r_2 & -r_3 & r_4 \end{bmatrix} \quad (11)$$

$$Q_{(r)}=\begin{bmatrix} r_4 I + C(r) & r \\ -r^T & r_4 \end{bmatrix} = \begin{bmatrix} r_4 & -r_3 & r_2 & r_1 \\ r_3 & r_4 & -r_1 & r_2 \\ -r_2 & r_1 & r_4 & r_3 \\ -r_1 & -r_2 & -r_3 & r_4 \end{bmatrix} \quad (12)$$

$W_{(r)}$ and $Q_{(r)}$ matrices are used to obtain the R transformation matrix from quaternions. For more detailed information about quaternions, readers are advised to look at the following studies Hamilton (1853), Wang et al. (2014), Zeng et al. (2020), Zeng et al. (2024). If the transformation will be made from unit quaternions, the following two equations with unity and orthogonality constraints must be added to the mathematical model.

$r^T r = 1$ , $r^T s = 0$ (13)

## 2.1 Symmetric DQA of 3D transformation model with constraint equations

Performing 3D coordinate transformations according to the dual quaternion and Euler angle methods are similar to each other. The functional model of 3D coordinate transformations based on Euler angle method is shown in Eq.(14).

$$X + v_{XYZ} = t + \lambda R (x + v_{xyz}) \quad (14)$$

The most important difference between them is that in the DQA method, the rotation matrix **R** and the translation vector **t** are obtained from the **r** and **s**. The functional model of symmetric DQA transformation equation is shown below Eq.(15). This model is used both for the transformation process and for calculating the coordinates of new points in the second system.

$$X + v_{XYZ} = 2W_{(r)}^T s + \lambda W_{(r)}^T Q_{(r)}(x + v_{xyz}) \quad (15)$$

Explicit form of the Eq. (15) with matrix representation

$$\begin{bmatrix} X + v_X \\ Y + v_Y \\ Z + v_Z \\ 0 \end{bmatrix}_{P_i} = 2W_{(r)}^T s + \lambda W_{(r)}^T Q_{(r)} \begin{bmatrix} x + v_x \\ y + v_y \\ z + v_z \\ 0 \end{bmatrix}_{P_i} \quad (16)$$

Transformation operations are performed on unit quaternions because this model contains a scale factor.

Where

$x$ vector denotes the coordinates of a point $P_i$ are $(x, y, z)$ in the first (source) system

$X$ vector denotes the coordinates of a point $P_i$ are $(X, Y, Z)$ in the second (target) system

$v_{xyz}$, $v_{XYZ}$ vector denotes the residual of $(x, y, z)$ and the residual of $(X, Y, Z)$ respectively

$\lambda$  scale factor between $(x, y, z)$ and $(X, Y, Z)$ coordinate systems

$t = 2W_{(r)}^T s$   translation vector

$$W_{(r)}^T Q_{(r)} = \begin{bmatrix} R_{(3\times3)} & 0_{(3\times1)} \\ 0_{(1\times3)} & 1 \end{bmatrix} \tag{17}$$

$W_{(r)}^T Q_{(r)}$ product matrix is the size of the $R_{(3\times3)}$ rotation matrix increased to $(4\times4)$ to be compatible with model Eq.(15-16).

Normally, four equations should be written per control point from this representation Eq. (15-16). When working with unit quaternions, an additional unity condition equation must be added Eq.(20). The total number of equations is *4n+1*. It should also be noted that the fourth equation ($r^T s = 0$) is written in the same way at each point. This equation does not carry any information about the point to which it belongs. For this reason, the fourth equation will not be written for each point. This equation, which is one of the orthogonality conditions, is written only once. As a result, the number of equations required decreases from *4n+1* to *3n+2*. The number of unknown parameters is nine in the constrained method.

If we write the transformation equation per point as an implicit observation equation,

$$f = X_i + v_{iXYZ} - 2W_{(r)}^T s - \lambda W_{(r)}^T Q_{(r)}(x_i + v_{ixyz}) \tag{18}$$

And two additional constraints for unity and orthogonality ($r^T s = 0$, $r^T r = 1$ )

$$r_1 s_1 + r_2 s_2 + r_3 s_3 + r_4 s_4 = 0 \tag{19}$$

$$r_1^2 + r_2^2 + r_3^2 + r_4^2 = 1 \tag{20}$$

$$r = [\, r_1 \; r_2 \; r_3 \; r_4 \,]^T \qquad s = [\, s_1 \; s_2 \; s_3 \; s_4 \,]^T \tag{21}$$

If the number of control points is three or more, adjustment is required. We cannot solve the symmetric transformation model by simple LS principle e.g. Gauss-Markov model. We need to use the below Error In Variables (EIV) model below.

$$A\, v + B\, x + w = 0 \tag{22}$$

However, two additional constraint equations need to be used. So our model will be Constraint Error In Variables (CEIV) Öztürk and Şerbetçi (1982).

$$A\, v + B\, x + w_1 = 0 \tag{23}$$

$$C\, x + w_2 = 0 \tag{24}$$

The transformation equation and constraint equation are nonlinear. It needs to be expanded to Taylor series and linearized concerning one scale factor and eight quaternions ( $\lambda$, $r_1$, $r_2$, $r_3$, $r_4$, $s_1$, $s_2$, $s_3$, $s_4$ ) also residuals ( $v_{x1}$ $v_{y1}$ $v_{z1}$ ...... $v_{xn}$ $v_{yn}$ $v_{zn}$ $v_{X1}$ $v_{Y1}$ $v_{Z1}$ ...... $v_{Xn}$ $v_{Yn}$ $v_{Zn}$).

We suppose that $\lambda_0$, $r_{1,0}$, $r_{2,0}$, $r_{3,0}$, $r_{4,0}$, $s_{1,0}$, $s_{2,0}$, $s_{3,0}$, $s_{4,0}$ are the approximate values of scale factor and quaternions.

Generally, the DQA method does not require appropriate approximations of the parameters for iterative computation. The approximations required for linearization can be obtained as follows.

$\lambda_0 = 1$, $r_{4,0} = 1$, $r_{1,0} = r_{2,0} = r_{3,0} = s_{1,0} = s_{2,0} = s_{3,0} = s_{4,0} = 0$

In exceptional cases, especially when the scale factor between the two systems is large, appropriate approximation may be required. Zeng et al. (2019) proposed an algorithm for appropriate approximation of *r* quaternions and $\lambda$ scale factor.

It was observed that the normal equations established in the symmetric transformation with the CEIV model are unstable and in extremely bad condition (of the order of $10^{-27}$). When trying to determine nine unknown parameters, it was seen that the scale factor and *r* quaternion converged easily, but the *s* quaternions did not converge and almost oscillated, meaning that the solution could not be obtained. This problem was also tried in networks with different structures, but it was observed that the problem continued to exist. Zeng et al. (2024) divided the unknowns into two groups for this problem. However, we did not see the need to divide the unknowns in our study. A different strategy was used to overcome this convergence problem. Instead of using the original coordinates in each iteration, we used the corrected coordinates with the residuals found in the previous iteration. This means changing the original measurements and is a dangerous situation. To eliminate this drawback, the contribution of the residuals added to the functional model at the beginning is eliminated by subtracting them from the misclosure vector. In other words, the convergence problem was solved by using the modified misclosure vector Eq.(47) Bektas (2024b), Mikhail and Ackermann (1976). The use of a modified misclosure vector is important. If it is not used, the iteration may stop at a different local minimum instead of the global minimum. The user may not realise that the results found are not the desired values. After linearization, unknown parameters.

$\delta \mathbf{x} = [\delta\lambda \ \delta r_1 \ \delta r_2 \ \delta r_3 \ \delta r_4 \ \delta s_1 \ \delta s_2 \ \delta s_3 \ \delta s_4]^T$ (25)

functional model

$\mathbf{A} \mathbf{v} + \mathbf{B} \delta \mathbf{x} + \mathbf{w}_1 = 0$ (26)

$\mathbf{C} \delta \mathbf{x} + \mathbf{w}_2 = 0$ (27)

Where

*n* : number of control points

$\mathbf{A}_{(3n \times 6n)}$ coefficient matrix of residuals

$\mathbf{B}_{(3n \times 9)}$ design matrix of $\delta \mathbf{x}$ parameter

$\mathbf{C}_{(2 \times 9)}$ coefficient matrix of constraint for $\delta \mathbf{x}$

$\mathbf{v}_{(6n \times 1)}$ is the residual vector both of two system

$v = [\ v_{x1}\ v_{y1}\ v_{z1}\ ......\ v_{xn}\ v_{yn}\ v_{zn}\ v_{X1}\ v_{Y1}\ v_{Z1}\ ......\ v_{Xn}\ v_{Yn}\ v_{Zn}]^T$

$\mathbf{w_1}_{(3n \times 1)}$ modified misclosure vector of residual equations

$\mathbf{w_2}_{(2 \times 1)}$ misclosure vector of constraint equations

$\mathbf{P_{xyz}}_{(3n \times 3n)}$, $\mathbf{P_{XYZ}}_{(3n \times 3n)}$ weight matrix of first $xyz$ system and second $XYZ$ system

The weight matrix is calculated from the variance-covariance $(C_{xyz})$ matrix of the measurements. Where $\sigma_1^2, \sigma_2^2$ represents the variance of the unit-weighted measurement of first $xyz$ system and second $XYZ$ system.

$$\mathbf{P_{xyz}} = \sigma_1^2\ C_{xyz}^{-1} \qquad\qquad \mathbf{P_{XYZ}} = \sigma_2^2\ C_{XYZ}^{-1} \tag{28}$$

$\mathbf{P}_{(6n \times 6n)}$ unified weight matrix

$$\mathbf{P} = \begin{bmatrix} P_{xyz} & 0 \\ 0 & P_{XYZ} \end{bmatrix} \tag{29}$$

The objective function for least squares adjustment

$$\Omega = v^T P v = \min. \tag{30}$$

In addition to minimizing the objective function, Eq. (26-27) in the model must also be satisfied. For this, we can add two Lagrange extreme conditions ($k_a$, $k_b$) that have zero effect on the objective function.

$$\Omega = v^T P v - 2\ k_a^T(\mathbf{A}\ v + \mathbf{B}\ \delta x + w_1) - 2\ k_b^T(\mathbf{C}\ \delta x + w_2) \tag{31}$$

For this function to be a minimum

$\frac{\partial \Omega}{\partial v}, \frac{\partial \Omega}{\partial k_a}, \frac{\partial \Omega}{\partial k_b}, \frac{\partial \Omega}{\partial \delta x}$   Its derivatives must be set equal to zero.

From the above derivatives,

$$\frac{\partial \Omega}{\partial v} = 2 v^T P - 2\ k_a^T \mathbf{A} = 0 \tag{32}$$

one get

$$v = P^{-1} A^T k_a \tag{33}$$

$$\frac{\partial \Omega}{\partial k_a} = A\ P^{-1} A^T k_a + B\ \delta x + w_1 = 0 \tag{34}$$

And introducing

$$N = A\ P^{-1} A^T \tag{35}$$

$$\frac{\partial \Omega}{\partial k_a} = N\ k_a + B\ \delta x + w_1 = 0 \tag{36}$$

$$\frac{\partial \Omega}{\partial k_b} = C\ \delta x + w_2 = 0 \tag{37}$$

$$\frac{\partial \Omega}{\partial \delta x} = B^T\ k_a + C^T\ k_b = 0 \tag{38}$$

and Lagrange multiplier's vector

$$k_{a(3n \times 1)} = [\ k_{a1}\ k_{a2} \ldots\ k_{a3n}\ ]^T \tag{39}$$

$$k_{b(2 \times 1)} = [\ k_{b1}\ \ k_{b2}\ ]^T \tag{40}$$

$$A_{(3n \times 6n)} = \begin{bmatrix} A_{1\,(3\times 3)} & & & & \vdots & \\ & A_{2(3\times 3)} & & & \vdots & \\ & & \vdots & & \vdots & I_{(3n \times 3n)} \\ & & & & \vdots & \\ & & & A_{n(3\times 3)} & \vdots & \end{bmatrix} \tag{41}$$

Asymmetric transformation results can also be obtained with a simple modification of the model established for symmetric transformation. For this, it is sufficient to simply change Eq. (41) as follows.

$$A_{(3n \times 6n)} = [\ 0_{(3n \times 3n)} \vdots\ I_{(3n \times 3n)}\ ] \tag{42}$$

It is possible to obtain asymmetric transformation results in other types of transformations (unconstrained DQA, QA) by using Eq.(42).

Where $0_{(3n \times 3n)}$ is a zero matrix and $I_{(3n \times 3n)}$ is an identity matrix.

$$B_{(3n \times 9)} = \begin{bmatrix} B_{1(3\times 9)} \\ B_{2(3\times 9)} \\ B_{3(3\times 9)} \\ \vdots \\ B_{n(3\times 9)} \end{bmatrix} \tag{43}$$

Submatrices $A_i$ and $B_i$ are obtained from the derivatives of the functional model Eq.(18)

$$A_{i(3\times 3)} = \left[ \frac{\partial f_i}{\partial(v_{xi},\ v_{yi},\ v_{zi})} \right] \tag{44}$$

$$B_{i(3\times9)} = \left[\frac{\partial f_i}{\partial(\lambda,r_1,r_2,r_3,r_4,s_1,s_2,s_3,s_4)}\right] \quad (45)$$

The desired derivatives in question can be easily achieved using the symbolic derivative functions of software such as Matlab.

The linearized constraint equations from Eq. (27)

$$C_{(2\times9)} = \begin{bmatrix} 0 & r_1 & r_2 & r_3 & r_4 & 0 & 0 & 0 & 0 \\ 0 & s_1 & s_2 & s_3 & s_4 & r_1 & r_2 & r_3 & r_4 \end{bmatrix} \quad (46)$$

$f$ misclosure vector from Eq. (18) and modified misclosure vector $w_1$

$$w_1^{(i)} = f^{(i)} - A^{(i)} \cdot v^{(i-1)} \quad (47)$$

$w_2$ misclosure vector of constraint equations

$$w_{2(2\times1)} = \begin{bmatrix} \frac{1}{2}(r_1^2+r_2^2+r_3^2+r_4^2 - 1) \\ r_1 s_1 + r_2 s_2 + r_3 s_3 + r_4 s_4 \end{bmatrix} \quad (48)$$

It should be noted that the original coordinates of both systems ($X$, $x$) must be corrected and used in each iteration.

$$\hat{X}_i = X + v^{(i-1)} \quad \hat{x}_i = x + v^{(i-1)} \quad (49)$$

$i$ denotes iteration number, for the first iteration, vector $v^{(0)} = [0\ 0\ 0...\ 0]^T$

If we rearrange these three equations Eqs.(36-38) and write

$$N\ k_a + 0\ k_b + B\ \delta x = -w_1 \quad (50)$$

$$0\ k_a + 0\ k_b + C\ \delta x = -w_2 \quad (51)$$

$$-B^T\ k_a - C^T\ k_b + 0\ \delta x = 0 \quad (52)$$

We express it in the matrix representation

$$M = \begin{bmatrix} N & 0 & B \\ 0 & 0 & C \\ -B^T & -C^T & 0 \end{bmatrix} \quad x = \begin{bmatrix} k_a \\ k_b \\ \delta x \end{bmatrix} \quad w = \begin{bmatrix} -w_1 \\ -w_2 \\ 0 \end{bmatrix} \quad (53)$$

$M\ x = w$ and $x = M^{-1}\ w \quad (54)$

The number of unknown parameters to be determined is $3n+11$.

$$x_{(3n+11)} = [\ k_{a1}\ k_{a2}\ ....\ k_{a3n}\ k_{b1}\ k_{b2}\ \delta\lambda\ \delta r_1\ \delta r_2\ \delta r_3\ \delta r_4\ \delta s_1\ \delta s_2\ \delta s_3\ \delta s_4\ ]^T \quad (55)$$

Using Eq. (54) $x$ can be obtained simultaneously in each iteration. The value of the scale factor and eight quaternions is improved by adding the calculated differential values by Eq. (56-58) at each iteration. The subscript ($i$) denotes present iterative times

$$r_{1.(i)} = r_{1.(i-1)} + \delta r_1 \quad r_{2.(i)} = r_{2.(i-1)} + \delta r_2 \quad r_{3.(i)} = r_{3.(i-1)} + \delta r_3 \quad r_{4.(i)} = r_{4.(i-1)} + \delta r_4 \tag{56}$$

$$s_{1.(i)} = s_{1.(i-1)} + \delta s_1 \quad s_{2.(i)} = s_{2.(i-1)} + \delta s_2 \quad s_{3.(i)} = s_{3.(i-1)} + \delta s_3 \quad s_{4.(i)} = s_{4.(i-1)} + \delta s_4 \tag{57}$$

$$\lambda_{(i)} = \lambda_{(i-1)} + \delta\lambda \tag{58}$$

Generally, the stop condition of iteration is $(\delta r^T \delta r + \delta s^T \delta s) < 10^{-11}$ when this condition is fulfilled, the iteration is terminated. Where $\delta r = [\delta r_1 \; \delta r_2 \; \delta r_3 \; \delta r_4]^T$ and $\delta s = [\delta s_1 \; \delta s_2 \; \delta s_3 \; \delta s_4]^T$

If desired scaled quaternion ($q_s$) calculated from unit quaternion ($r$)

$$r_1^2 + r_2^2 + r_3^2 + r_4^2 = 1 \tag{59}$$

$$q_{si} = r_i \sqrt{\lambda} \qquad (i=1,2,3,4) \tag{60}$$

The residuals to the control points are calculated via Eq. (61).

$$v = P^{-1} A^T k_a \tag{61}$$

A posteriori standard error of an observation of unit weight is given by

$$\widehat{\sigma_0} = \sqrt{\frac{v^T P v}{f}} = \sqrt{\frac{v^T P v}{3n-7}} \tag{62}$$

Where $f$ denotes the degree of freedom of transformation.

The flow of the presented similarity symmetric 3D coordinate transformation based on a dual quaternion algorithm is summarized in Table 1.

**Table 1 Similarity symmetric transformation based on dual quaternion algorithm**

---

**Input and initiation:**

input 3D coordinates of control points (**xyz**) and (**XYZ**), weight matrix or variance of points initial value $\lambda_0=1$, $r_{4,0}=1$, $r_{1,0} = r_{2,0} = r_{3,0} = s_{1,0} = s_{2,0} = s_{3,0} = s_{4,0} = 0$ and $v^{(0)} =[0\ 0\ ....0]^T$

**Step 1** Construct corrected both of coordinates system

$\hat{X}_i = X + v^{(i-1)} \qquad \hat{x}_i = x + v^{(i-1)}$

The subscript *i* denotes present iterative times. Construct mathematic model by Eq. (26-27)

Calculate **P, A, B , C , w₁ ,w₂** by Eq.(29) and Eqs. (44-48).

Construct Normal Equations $M\ x = w$ by Eq.(54)

**Step 2** Compute iteratively unknown parameters by Eq.(54) ($k_{a1}\ k_{a2}....k_{a3n}\ k_{b1}\ k_{b2}\ \delta\lambda\ \delta r_1\ \delta r_2\ \delta r_3\ \delta r_4\ \delta s_1\ \delta s_2\ \delta s_3\ \delta s_4$)

**Step 3** Improve approximations

$r_i = r_{i,0} + \delta r_i$ , $s_i = s_{i,0} + \delta s_i$ , $\lambda_i = \lambda_0 + \delta\lambda_i$ by Eq. (56-58)

**Step 4** If $(\delta r^T \delta r + \delta s^T \delta s) < 10^{-11}$ turn to Step 1 otherwise turn to Step 5

**Step 5** Compute residual by Eq.(58) then a posteriori standard deviation $\widehat{\sigma_0} = \sqrt{\frac{v^T P v}{3n-7}}$ by Eq.(62).

**Step 6** If desired compute scaled quaternions ($q_s$) by Eq. (60). Compute **R** rotation matrix substituting *r* quaternions instead of *q* quaternions in Eq. (8)

**Step 7** Check the results; convert the coordinates using the calculated transformation parameters from first system to the second system. By adding corrections to the coordinates of both systems. Make sure that the corrected coordinates of the second system are exactly obtained from the corrected coordinates of first system by the transformation equation of Eq. (15).

**Step 8** If there are new points to be converted, compute new points coordinate in the second system by Eq. (69-70).

**Step 9** If rotation angles (**ε, ψ, ω**) and translations parameters($t_X$ , $t_Y$, $t_Z$) are needed by Eq. (75-81), compute the variance of all the transformation parameters ( $r_1\ r_2\ r_3\ r_4\ s_1\ s_2\ s_3\ s_4$ ) by Eq.(71-74), and compute the variance of six transformation parameters ($\lambda, \varepsilon, \psi, \omega,\ t_X, t_Y, t_Z$) .

**Using scaled quaternions**

The transformation operations we carry out on unit quaternions can also be carried out on scaled quaternions in a similar way. In this case, what needs to be done can be summarized as follows. The scale parameter is removed from the transformation equation Eq.(18).

$$f_i = X_i + v_{iXYZ} - 2W_{(r)}^T s - W_{(r)}^T Q_{(r)}(x_i + v_{ixyz}) \tag{63}$$

In this case, the number of additional constraint to be written is only one ($r^T s = 0$). The total number of equations to be written is $3n+1$. The number of unknown parameters becomes eight ($r_1\ r_2\ r_3 r_4\ s_1 s_2\ s_3\ s_4$). Since the scale factor is removed from the unknowns parameters $x$, $B$ and $C$ matrix have eight columns. In Eq.(46) and Eq.(48) the first lines of $C$ and $w_2$ are deleted as below.

$$\mathbf{x} = [r_1\ r_2\ r_3\ r_4 s_1 s_2\ s_3\ s_4]^T \tag{64}$$

$$\mathbf{B}_{i(3\times8)} = \left[\frac{\partial f_i}{\partial(r_1,r_2,r_3,r_4,s_1,s_2,s_3,s_4)}\right] \tag{65}$$

$$\mathbf{C}_{(1\times8)} = [s_1\quad s_2\quad s_3\quad s_4\quad r_1\quad r_2\quad r_3\quad r_4] \tag{66}$$

$$\mathbf{w}_{2(1\times1)} = [r_1 s_1 + r_2 s_2 + r_3 s_3 + r_4 s_4] \tag{67}$$

After calculating the scaled quaternions, the scale factor is calculated as follows.

$$\lambda = r_1^2 + r_2^2 + r_3^2 + r_4^2 \tag{67}$$

The iterative calculation and other precision calculations are carried out as in Table 1.

**Checking results**

It is a general rule in surveying that measurements and calculations are made in a controlled way. The cost of making a mistake in professional life can be very high. For this reason, both the measurement and the calculation processes should be checked. When the iterative calculation is finished, both system coordinates are corrected by calculating the corrections (residuals), and the corrected coordinates must fully satisfy the transformation equation Eq. (69).

$$X + v_{XYZ} \stackrel{?}{=} 2W_{(r)}^T s + \lambda W_{(r)}^T Q_{(r)}(x + v_{xyz}) \tag{69}$$

If not, the reason should be investigated. If there are new points to be converted, compute new points coordinate in the second system by Eq. (73)

$$\begin{bmatrix} X \\ Y \\ Z \\ 0 \end{bmatrix}_{new} = 2W_{(r)}^T s + \lambda\, W_{(r)}^T Q_{(r)} \begin{bmatrix} x \\ y \\ z \\ 0 \end{bmatrix}_{new} \qquad (70)$$

**2.2. Computation and precision estimation of six parameters from dual quaternion**

The precision values for *3n+11* unknown parameters (*3n+2* Lagrange multipliers + *1* scale factor + *8* quaternions) can be easily calculated directly by using the inverse of the **M** Matrix Eq. (56). Although one scale factor and eight quaternions are sufficient to perform the symmetric similarity transformation calculations, it is still desirable to determine the classical six geometrical transformation parameters (three translations, three rotation angles) due to habits. In this section, we will see how to determine six parameters from quaternions and how to calculate the precision of the calculated parameters. Using the inverse of the (**M**$^{-1}$) matrix, the full precision information of all parameters can be calculated. Since we are using unit quaternions, precision calculations are easy to do. A posteriori standard variance $\sigma_0^2$ from Eq.(62)

$$C = \sigma_0^2\, M^{-1} = \begin{bmatrix} C_{k(3n+2)x(3n+2)} & : \\ : & C_{x(9\times 9)} \end{bmatrix} \qquad (71)$$

$C_{(3n+11\times 3n+11)}$ is the variance covariance matrix of all parameters including (3n+2) Lagrange multipliers, one scale factor and eight quaternions.

The standard error of the scale factor is

$$\sigma_\lambda = \sqrt{C_{3n+3,3n+3}} \qquad (72)$$

The standard error of the *r* quaternions are

$\sigma_{ri}=\sqrt{C_{j,j}}$     ($i$=1,2,3,4)  ($j$=$i$+3n+3)     (73)

The standard error of the *s* quaternions are

$\sigma_{si}=\sqrt{C_{j,j}}$     ($i$=1,2,3,4)  ($j$=$i$+3n+7)     (74)

To calculate the precision of six transformation parameters, first, the relational equations between six transformation parameters and eight quaternions are written as follows. As the scale parameter is determined directly in the model, the number of parameters to be determined is only six (three rotation angles, three translations).

Euler rotation angles can be calculated from the quaternions using the following equations.

*ε = -atan2 (2($r_4 r_1$+$r_2 r_3$) , ($r_4^2$-$r_1^2$-$r_2^2$+$r_3^2$))*     (75)

*ψ = arcsin (2 ($r_3 r_1$-$r_4 r_2$))*     (76)

*ω = -atan2 (2 ($r_4 r_3$+$r_2 r_1$) , ($r_4^2$+$r_1^2$-$r_2^2$-$r_3^2$))*     (77)

The translations parameters

$t = 2W_{(r)}^T s = [t_X\ t_Y\ t_Z\ 0]^T$     (78)

$t_X = 2(r_2 s_3 - r_1 s_4 - r_3 s_2 + r_4 s_1)$     (79)
$t_Y = 2(r_3 s_1 - r_1 s_3 - r_2 s_4 + r_4 s_2)$     (80)
$t_Z = 2(r_1 s_2 - r_2 s_1 - r_3 s_4 + r_4 s_3)$     (81)

The variance-covariance matrix of the six parameters according to the law of propagation of variances

$cx = \cos\varepsilon$     $cy = \cos\psi$     $cz = \cos\omega$     (82)
$sx = \sin\varepsilon$     $sy = \sin\psi$     $sz = \sin\omega$     (83)

$$J_\alpha = \left[\frac{\partial(\varepsilon,\psi,\omega)}{\partial(r_1,r_2,r_3,r_4)}\right] = \frac{-1}{cy}\begin{bmatrix} -r_1.cy & -r_2.cy & -r_3.cy & -r_4.cy \\ (r_4.cx - r_1.sx) & (r_3.cx - r_2.sx) & (r_2.cx + r_3.sx) & (r_1.cx + r_4.sx) \\ (-r_3 + r_1.sy) & (r_4 + r_2.sy) & (-r_1 + r_3.sy) & (r_2 + r_4.sy) \\ (r_2.cz + r_1.sz) & (r_1.cz - r_2.sz) & (r_4.cz - r_3.sz) & (r_3.cz + r_4.sz) \end{bmatrix}$$ (84)

Uygur et al. (2022)

$$J_t = \left[\frac{\partial(t_X,t_Y,t_Z)}{\partial(r_1,r_2,r_3,r_4,s_1,s_2,s_3,s_4)}\right] = 2\begin{bmatrix} -s_4 & s_3 & -s_2 & s_1 & r_4 & -r_3 & r_2 & -r_1 \\ -s_3 & -s_4 & s_1 & s_2 & r_3 & r_4 & -r_1 & -r_2 \\ s_2 & -s_1 & -s_4 & s_3 & -r_2 & r_1 & r_4 & -r_3 \end{bmatrix}$$ (85)

Adopted Bektas (2024a)

$$J = \begin{bmatrix} J_\alpha & 0 \\ & J_t \end{bmatrix} \quad (86)$$

$J_{\alpha(4\times 4)}$, $0_{(4\times 4)}$, $J_{t(3\times 8)}$ and $J_{(7\times 8)}$ are Jacobian matrix Bektas (2024a).

$$C_{(3n+11\times 3n+11)} = \begin{bmatrix} C_{kk(3n+3\times 3n+3)} & C_{kq(3n+3\times 8)} \\ C_{qk(8\times 3n+3)} & C_{qq(8\times 8)} \end{bmatrix} \quad (87)$$

Let matrix $C_{qq}$ be the last 8×8 sub diagonal matrix of matrix $C$. By applying the variances-covariance matrix propagation principle. The variance-covariance matrix of these six transformation parameters $C_{par}$

$$C_{par} = J \, C_{qq} \, J^T \quad (88)$$

The diagonal elements of the $C_{par}$ matrix are the variances of the parameters, respectively $\varepsilon, \psi, \omega, t_X, t_Y, t_Z$

The presented 3D symmetric similarity coordinate transformation based on a dual quaternion algorithm is summarized in Table 1

## 2.3. Simplify symmetric DQA of 3D transformation model

Symmetric 3D coordinate transformations can also be performed with the classical EIV model without the need for constraint equations. For example, for the $r^T r = 1$ unity constraint, the $r_4$ quaternion is removed, and for the $r^T s = 0$ orthogonality constraint, the $s_4$ quaternion is removed from the unknown parameters in the functional model. Instead of the removed $r_4$, $s_4$ quaternions

$$r_4 = (1 - r_1^2 - r_2^2 - r_3^2)^{1/2} \quad (89)$$

$$s_4 = -(r_1 s_1 + r_2 s_2 + r_3 s_3) / (1 - r_1^2 - r_2^2 - r_3^2)^{1/2} \quad (90)$$

is written in the functional model.

Linearized functional model from Eq. (22).

$$f = A_{(3n\times 6n)} \, v_{(6n\times 1)} + B_{(3n\times 7)} \, \delta x_{(7\times 1)} + w_{(3n\times 1)} = 0_{(3n\times 1)} \quad (91)$$

$\delta x_{(7\times1)} = [\ \delta r_1\ \delta r_2\ \delta r_3\ \delta s_1\ \delta s_2\ \delta s_3\ \lambda\ ]^T$  independent unknowns parameter

**A** and **B** matrices are created in accordance with Eq.(41-43)

$$\mathbf{A}_{(3n\times 6n)} = \begin{bmatrix} \mathbf{A}_{1\,(3\times 3)} & & & & : & \\ & \mathbf{A}_{2(3\times 3)} & & & : & \\ & & : & & : & \mathbf{I}_{(3n\times 3n)} \\ & & & & : & \\ & & & \mathbf{A}_{n(3\times 3)} & : & \end{bmatrix} \quad (92)$$

$$\mathbf{B}_{(3n\times 7)} = \begin{bmatrix} \mathbf{B}_{1(3\times 7)} \\ \mathbf{B}_{2(3\times 7)} \\ \mathbf{B}_{3(3\times 7)} \\ : \\ \mathbf{B}_{n(3\times 7)} \end{bmatrix} \quad (93)$$

The modified misclosure vector computed as Eq. (47)

$w_1^{(i)} = f^{(i)} - A^{(i)} \cdot v^{(i-1)}$   (The subscript $i$ denotes present iterative times)

Submatrices of **A** and **B** matrices

$A_{i(3\times 3)} = \left[\dfrac{\partial f_i}{\partial(v_{xi},\ v_{yi},\ v_{zi})}\right] \qquad B_{i(3\times 7)} = \left[\dfrac{\partial f_i}{\partial(r_1,r_2,r_3,s_1,s_2,s_3,\lambda)}\right]$   (i=1,2,..,n)   **(94)**

$r_4$, $s_4$ the dependent unknowns parameter computed from Eq.(89-90)

Stochastic model from Eq. (29). The objective function for least squares adjustment

$\Omega = v^T P v = $ min.

In addition to minimizing the objective function, Eq. (91) in the model must also be satisfied. For this, we can add one Lagrange extreme conditions (***k***) that have zero effect on the objective function.

$\Omega = v^T P v - 2\ k^T(A\ v + B\ x + w) = \boldsymbol{min.}$   **(95)**

For $\Omega$ function to be a minimum $\dfrac{\partial \Omega}{\partial v},\dfrac{\partial \Omega}{\partial k},\dfrac{\partial \Omega}{\partial \delta x}$   Its derivatives must be set equal to zero.

$N = A\ P^{-1} A^T$   **(96)**

unknown parameters

$\delta x = -\ (B^T N^{-1} B)^{-1}\ B^T N^{-1} w$   **(97)**

Lagrange multiplier

$$k = - N^{-1} (B \, \delta x + w) \tag{98}$$

The residuals to the control points are calculated via Eq. (99).

$$v = P^{-1} A^T k \tag{99}$$

A posteriori standard error of an observation of unit weight is given by

$$\widehat{\sigma_0} = \sqrt{\frac{v^T P v}{f}} = \sqrt{\frac{v^T P v}{3n-7}} \tag{100}$$

The variance-covariance matrix of all parameters ($r_1$, $r_2$, $r_3$, $s_1$, $s_2$, $s_3$, $\lambda$) six quaternions and one scale factor.

$$C_{(7 \times 7)} = \sigma_0^2 \, (B^T N^{-1} B)^{-1} \tag{101}$$

The error of the dependent parameter ($r_4$, $s_4$) is found by applying the law of propagation of variances to Eq.(89-90).

$$J_{r4s4} = \begin{bmatrix} \partial r_4/\partial r_1 & \partial r_4/\partial r_2 & \partial r_4/\partial r_3 & \partial r_4/\partial s_1 & \partial r_4/\partial s_2 & \partial r_4/\partial s_3 & 0 \\ \partial s_4/\partial r_1 & \partial s_4/\partial r_2 & \partial s_4/\partial r_3 & \partial s_4/\partial s_1 & \partial s_4/\partial s_2 & \partial s_4/\partial s_3 & 0 \end{bmatrix} \tag{102}$$

$$C_{r4s4} = \begin{bmatrix} C_{r4r4} & C_{r4s4} \\ C_{s4r4} & C_{s4s4} \end{bmatrix} = J_{r4s4} \, C_{(7 \times 7)} \, J^T_{r4s4} \tag{103}$$

Error of $r$ quaternions

$$\sigma_{ri} = \sqrt{C_{i,i}} \quad (i=1,2,3) \qquad \sigma_{r4} = \sqrt{C_{r4r4}} \tag{104}$$

Error of $s$ quaternions

$$\sigma_{si} = \sqrt{C_{j,j}} \quad (i=1,2,3) \quad (j=i+3) \qquad \sigma_{s4} = \sqrt{C_{s4s4}} \tag{105}$$

Error of $\lambda$ scale factor

$$\sigma_\lambda = \sqrt{C_{7,7}} \tag{106}$$

Euler rotation angles and translation parameters can be calculated from quaternions using the Eq. (75-81). To calculate the precision of the six transformation parameters $(\varepsilon, \psi, \omega, t_X, t_Y, t_Z)$, firstly, the Jacobian matrix ($J_p$) is calculated by taking derivatives according to the six quaternions over the equations Eq.(75-77) and Eq.(79-81). When taking derivatives, it should not be forgotten that instead of $(r_4, s_4)$ quaternions, their equivalents in terms of others are written as in Eq.(89-90).

$$J_p = \left[\frac{\partial(\varepsilon, \psi, \omega, t_X, t_Y, t_Z)}{\partial(r_1, r_2, r_3, s_1, s_2, s_3)}\right] = \begin{bmatrix} \partial\varepsilon/\partial r_1 & \partial\varepsilon/\partial r_2 & \partial\varepsilon/\partial r_3 & \partial\varepsilon/\partial s_1 & \partial\varepsilon/\partial s_2 & \partial\varepsilon/\partial s_3 \\ \partial\psi/\partial r_1 & \partial\psi/\partial r_2 & \partial\psi/\partial r_3 & \partial\psi/\partial s_1 & \partial\psi/\partial s_2 & \partial\psi/\partial s_3 \\ \partial\omega/\partial r_1 & \partial\omega/\partial r_2 & \partial\omega/\partial r_3 & \partial\omega/\partial s_1 & \partial\omega/\partial s_2 & \partial\omega/\partial s_3 \\ \partial t_X/\partial r_1 & \partial t_X/\partial r_2 & \partial t_X/\partial r_3 & \partial t_X/\partial s_1 & \partial t_X/\partial s_2 & \partial t_X/\partial s_3 \\ \partial t_Y/\partial r_1 & \partial t_Y/\partial r_2 & \partial t_Y/\partial r_3 & \partial t_Y/\partial s_1 & \partial t_Y/\partial s_2 & \partial t_Y/\partial s_3 \\ \partial t_Z/\partial r_1 & \partial t_Z/\partial r_2 & \partial t_Z/\partial r_3 & \partial t_Z/\partial s_1 & \partial t_Z/\partial s_2 & \partial t_Z/\partial s_3 \end{bmatrix} \quad (107)$$

by applying the law of propagation of variances to Eq. (79-81) and Eq. (83-85).

$$\boldsymbol{C_p} = \boldsymbol{J_p}\, C_{(6\times6)}\, \boldsymbol{J_p^T} \qquad (108)$$

$\boldsymbol{C_p}$ is the variance-covariance matrix of six transformation parameters. $C_{(6\times6)}$ matrix is the part of $C_{(7\times7)}$ corresponding to six parameters. In other words, $C_{(6\times6)}$ is obtained by deleting the last row and last column of $C_{(7\times7)}$.

Error of rotation angles

$\sigma_\varepsilon = \sqrt{C_{p1,1}} \qquad \sigma_\psi = \sqrt{C_{p2,2}} \qquad \sigma_\omega = \sqrt{C_{p3,3}} \qquad$ [in radian] $\qquad (109)$

Error of translation parameters

$\sigma_{tX} = \sqrt{C_{p4,4}} \qquad \sigma_{tY} = \sqrt{C_{p5,5}} \qquad \sigma_{tZ} = \sqrt{C_{p6,6}} \qquad$ [in meter] $\qquad (110)$

The iterative calculation and other precision calculations are carried out as in Table 1, the constraint algorithm previously presented.

This, the number of unknown parameters decreases from nine to seven (3 *r* quaternions + 3 *s* quaternions + 1 scale factor). The total number of equations is also reduced from *3n+2* to *3n*. The advantage of this simple method is that it reduces the number of unknown parameters and the total number of equations. There is also a simpler model that is easy to solve (there are no constraint equations).

## 2.4. Quaternion Algorithm (QA) of symmetric 3D transformation model with unit quaternions

Symmetric 3D coordinate transformations can also be performed with the QA model using unit quaternions. The purpose of adding this section is to make the contradictions in the precision calculations in Zeng's study clearer. The function model of the symmetric 3D transformation model with unit quaternions is below. By analogy from Eq. (18)

$$\boldsymbol{f}_i = \boldsymbol{X}_i + \boldsymbol{v}_{iXYZ} - [t_X\ t_Y\ t_Z\ 0]^T - \lambda \boldsymbol{W}_{(r)}^T \boldsymbol{Q}_{(r)}(\boldsymbol{x}_i + \boldsymbol{v}_{ixyz}) \tag{111}$$

The main difference between QA and DQA functional models appears to be in the translation vector Eq. 82). The approximate values of the unknown parameters required for linearization can be taken as follows. $\lambda_0 = 1$, $r_{4,0} = 1$, $r_{1,0} = r_{2,0} = r_{3,0} = t_X = t_Y = t_Z = 0$

Without using additional constraint equations, a functional model is created by simply substituting the following instead of $r_4$ as follows.

$r_4 = (1 - r_1{}^2 - r_2{}^2 - r_3{}^2)^{1/2}$

In this case, the number of independent unknowns is 7 (one scale factor + three translations + three *r* quaternions). The solution is continued as in section 5. The only difference from the functional model in Section 5 is in the derivative submatrix ($\boldsymbol{B}_i$) with respect to the unknowns.

$$B_{i(3\times 7)} = \left[\frac{\partial f_i}{\partial(t_X, t_Y, t_Z, r_1, r_2, r_3, \lambda)}\right] \qquad (i=1,2,..,n) \qquad (112)$$

The iterative calculation and other precision calculations are carried out in Table 1.

As can be seen from numerical applications, the QA and DQA transformation models give results that are fully compatible with each other.

## 3. Numerical experiments and discussion

Two numerical examples are designed to verify the correctness and effectiveness of the two DQA solutions presented. All examples are taken from Zeng et al. (2024). The results of the first two numerical examples were compared with the results of Uygur et al. (2020) and Zeng et al. (2024). To better explain the contradiction in the precision calculations of Zeng et al. (2024), the results of a new QA algorithm using unit quaternions and a DQA algorithm using scaled quaternions have also been added to the table. The transformation parameters, their standard errors and all the variance-covariance matrices are shown in the tables. All calculations were performed in the two algorithm model (with constraint and without constraint equation). At the same time, all results audits were carried out. It was observed that exactly the same results (transformation parameters and precision values) were achieved with the two methods we proposed.

### Case study 1

The data are chosen (Actual geodetic datum transformation case) by Zeng et al. (2024), Grafarend and Awange (2003). The 3D coordinates of the control point in (xyz) 1<sup>st</sup> and (XYZ) 2<sup>nd</sup> system and the point's variance are listed in Table 2

**Table 2** Coordinate of control points (xyz) and (XYZ) system and the point's variance

| x[m] | y[m] | z[m] | $\sigma^2_{xyz}$ | X[m] | Y[m] | Z[m] | $\sigma^2_{XYZ}$ |
|---|---|---|---|---|---|---|---|
| 4157222.5430 | 664789.3070 | 4774952.0990 | 0.1433 | 4157870.2370 | 664818.6780 | 4775416.5240 | 0.0103 |
| 4149043.3360 | 688836.4430 | 4778632.1880 | 0.1551 | 4149691.0490 | 688865.7850 | 4779096.5880 | 0.0038 |
| 4172803.5110 | 690340.0780 | 4758129.7010 | 0.1503 | 4173451.3540 | 690369.3750 | 4758594.0750 | 0.0006 |
| 4177148.3760 | 642997.6350 | 4760764.8000 | 0.1400 | 4177796.0640 | 643026.7000 | 4761228.8990 | 0.0114 |
| 4137012.1900 | 671808.0290 | 4791128.2150 | 0.1459 | 4137659.5490 | 671837.3370 | 4791592.5310 | 0.0068 |
| 4146292.7290 | 666952.8870 | 4783859.8560 | 0.1469 | 4146940.2280 | 666982.1510 | 4784324.0990 | 0.00002 |
| 4138759.9020 | 702670.7380 | 4785552.1960 | 0.1220 | 4139407.5060 | 702700.2270 | 4786016.6450 | 0.0041 |

**Table 3** The computed parameters and their precision by QA and DQA 3D transformation

| | Quaternion 3D transformation (QA) | | Dual quaternion 3D transformation (DQA) | | |
|---|---|---|---|---|---|
| | Uygur et al.2022 (Scaled quat.) (1) | Ours (Unit quat.) (2) | Zeng's study (Table4-5) (Unit quat.) (3) | Presented Methods (Scaled quat.) (4) | (Unit quat.) (5) |
| $t_X$[m] | 641.8395 ±9.0327 | 641.8395 ±9.0327 | 641.8395 **±7.91832** | 641.8395 ±9.0327 | 641.8395 ±9.0327 |
| $t_Y$[m] | 68.4729 ±10.5317 | 68.4729 ±10.5317 | 68.4728 **±10.50456** | 68.4729 ±10.5317 | 68.4729 ±10.5317 |
| $t_Z$[m] | 416.2156 ±9.0495 | 416.2156 ±9.0495 | 416.2155 **±7.50345** | 416.2156 ±9.0495 | 416.2156 ±9.0495 |
| $\lambda$ | 1.00000561109 ±0.00000108 | 1.00000561109 ±0.00000108 | 1.00000561109 **0.000000007** | 1.00000561109 ±0.00000108 | 1.00000561109 ±0.00000108 |
| $\varepsilon$ [°] | -0.000277143 ±0.00008517 | -0.000277143 ±0.00008517 | -0.000277143 ±0.00008517 | -0.000277143 ±0.00008517 | -0.000277143 ±0.00008517 |
| $\psi$ [°] | 0.000248913 ±0.00009628 | 0.000248913 ±0.00009629 | 0.000248913 ±0.00009628 | 0.000248913 ±0.00009629 | 0.000248913 ±0.00009629 |
| $\omega$ [°] | 0.000273856 ±0.00007552 | 0.000273857 ±0.00007552 | 0.000273857 ±0.00007552 | 0.000273857 ±0.00007552 | 0.000273857 ±0.00007552 |
| $r_4$ | 1.000002806 ±5.41461e-7 | 0.9999999999 ±3.913e-12 | 0.9999999999 ±3.913e-12 | 1.000002806 ±5.41461e-7 | 0.9999999999 ±3.913e-12 |
| $r_1$ | 0.00000241852 ±7.43266e-7 | 0.00000241852 ±7.4327e-7 | 0.000002419 ±7.4327e-7 | 0.000002419 ±7.4326463e-7 | 0.00000241852 ±7.43265e-7 |
| $r_2$ | -.00000217218 ±8.40281e-7 | -0.0000021722 ±8.40281e-7 | -0.000002172 ±8.49313e-7 | -0.000002172 ±8.40281e-7 | -0.0000021722 ±8.40281e-7 |
| $r_3$ | -.00000238984 ±6.590306e-7 | -.00000238984 ±6.590298e-7 | -0.0000023898 ±6.5898e-7 | -0.000002390 ±6.590298e-7 | -0.00000238984 ±6.590298e-7 |
| $s_4$ | | | 320.9201 **±3.9593** | 320.919239 ±4.5166 | 320.92019 ±4.5166 |
| $s_1$ | | | 34.2377 **±5.2526** | 34.237601 ±5.2662 | 34.2377 ±5.2662 |
| $s_2$ | | | 208.106983 **±3.7520** | 208.106434 ±4.5250 | 208.10703 ±4.5250 |
| $s_3$ | | | -0.00020443 ±0.00022 | -0.0002044 ±0.000218 | -0.00020440 ±0.000218 |
| $\widehat{\sigma_0}$[m] | 0.1976 | 0.1976 | 0.1976 | 0.1976 | 0.1976 |

In column (3), different values of Zeng's study are shown in bold in the table. It is considered that the differences may result from a calculation error or typos. It can be seen that the results of the QA-based 3D transformation we performed with unit quaternions column (2) are fully compatible with the DQA-based results column (5). In this example, the scale factor being very close to 1 brings the scaled and unit quaternion values very close to each other. The computed residual of coordinates in the first and second systems by QA and DQA methods are shown in Table 5. The variance-covariance matrix of all parameters is shown in Table 5-6.

The residuals are shown in Table 4.

**Table 4** Residuals

| Points | Residual[m] | | | | | |
|---|---|---|---|---|---|---|
| | $v_x$ | $v_y$ | $v_z$ | $v_X$ | $v_Y$ | $v_Z$ |
| Solitude | -0.0885 | -0.1261 | -0.1313 | 0.0064 | 0.0091 | 0.0094 |
| Buoch Zeil | -0.0593 | 0.0489 | -0.0140 | 0.0015 | -0.0012 | 0.0003 |
| Hohenneufen | 0.0386 | 0.0887 | 0.0071 | -0.0002 | -0.0004 | -0.0000 |
| Kuehlenberg | -0.0181 | 0.0203 | 0.0803 | 0.0015 | -0.0017 | -0.0065 |
| Ex Mergelaec | 0.0860 | -0.0138 | 0.0049 | -0.0040 | 0.0006 | -0.0002 |
| Ex Hof Asperg | 0.0105 | -0.0069 | 0.0542 | -0.0000 | 0.0000 | -0.0000 |
| Ex Kaisersbach | 0.0257 | -0.0035 | -0.0022 | -0.0009 | 0.0001 | 0.0001 |

**Table 5** The part of *C* variance-covariance matrix belonging to the ($\lambda$, $r_1$, $r_2$, $r_3$, $r_4$, $s_1$, $s_2$, $s_3$, $s_4$) parameters

```
 1.173e-12   2.067e-24  6.201e-24 -2.067e-24  1.577e-29 -2.436e-06 -3.965e-07 -2.801e-06 -1.664e-12
 3.821e-26   5.524e-13 -2.404e-13 -1.89e-13  -2.31e-18   1.021e-06  3.424e-06 -1.372e-06 -1.28e-10
-2.36e-25   -2.404e-13  7.061e-13  1.498e-13  2.473e-18 -3.271e-06 -1.771e-06  3.096e-06  3.326e-11
-1.561e-26  -1.89e-13   1.498e-13  4.343e-13  1.82e-18  -4.221e-07 -2.707e-06  7.502e-07 -3.793e-11
 1.654e-24  -2.31e-18   2.473e-18  1.82e-18   1.531e-23 -1.058e-11 -1.86e-11   1.184e-11  2.913e-16
-2.436e-06   1.021e-06 -3.271e-06 -4.221e-07 -1.058e-11  20.4       7.451     -8.462     -0.0001811
-3.965e-07   3.424e-06 -1.771e-06 -2.707e-06 -1.86e-11   7.451      27.73     -8.724     -0.0004535
-2.801e-06  -1.372e-06  3.096e-06  7.502e-07  1.184e-11 -8.462     -8.724     20.48      0.0002287
-1.664e-12  -1.28e-10   3.326e-11 -3.793e-11  2.91e-16  -0.0001811 -0.0004535  0.0002287 4.78e-08
```

**Table 6** The variance-covariance *Cpar* matrix ($\varepsilon, \psi, \omega, t_X, t_Y, t_Z$)

```
  2.21e-12  -9.617e-13 -7.559e-13 -4.082e-06 -1.369e-05  5.488e-06
 -9.617e-13  2.824e-12  5.994e-13  1.308e-05  7.083e-06 -1.238e-05
 -7.559e-13  5.994e-13  1.737e-12  1.688e-06  1.083e-05 -3.001e-06
 -4.082e-06  1.308e-05  1.688e-06  81.59      29.8      -33.84
 -1.369e-05  7.083e-06  1.083e-05  29.8      110.9      -34.89
  5.488e-06 -1.238e-05 -3.001e-06 -33.84     -34.89      81.89
```

## Case study 2

The data are chosen (Simulated case) by Zeng et al. (2024). The 3D coordinates of the control point in (xyz) 1st and (XYZ) 2nd system and the point's weight are listed in Table 7

**Table 7** Coordinate of control points (xyz) and (XYZ) system

| Point | x[m] | y[m] | z[m] | X[m] | Y[m] | Z[m] | P weight |
|---|---|---|---|---|---|---|---|
| 1 | 30 | 40 | 10 | 290 | 150 | 15 | 1 |
| 2 | 100 | 40 | 10 | 420 | 80 | 2 | 2 |
| 3 | 100 | 130 | 10 | 540 | 200 | 20 | 2.5 |
| 4 | 30 | 130 | 10 | 390 | 300 | 5 | 4 |

**Table 8** Computed parameters and their precision by QA and DQA 3D transformation method

| | Quaternion 3D transformation (QA) | | Dual quaternion 3D transformation (DQA) | | |
|---|---|---|---|---|---|
| | Uygur et al.2022 (Scaled quat.) (1) | Ours (Unit quat.) (2) | Zeng's study (Table 11) (3) | Presented methods (Scaled quat.) (4) | (Unit quat.) (5) |
| $t_X[m]$ | 192.2444 ±20.2709 | 192.2444 ±20.2709 | 192.2444 ±20.2709 | 192.2444 ±20.2709 | 192.2444 ±20.2709 |
| $t_Y[m]$ | 109.9534 ±20.1299 | 109.9534 ±20.1299 | 109.9534 ±20.1299 | 109.9534 ±20.1299 | 109.9534 ±20.1299 |
| $t_Z[m]$ | -24.0823 ±29.06571 | -24.0823 ±29.06571 | -24.0823 ±29.06571 | -24.0823 ±29.06571 | -24.0823 ±29.06571 |
| $\lambda$ | 2.136189318 ±0.152489951 | 2.1361893188 ±0.152489951 | 2.1361893188 ±0.152489951 | 2.1361893188 ±0.152489951 | 2.136189318 ±0.152489951 |
| $\varepsilon$ [°] | -1.882226178 ±5.8810538 | -1.882226178 ±5.8810538 | -1.882226178 ±5.8810538 | -1.882226178 ±5.8810538 | -1.882226178 ±5.8810538 |
| $\psi$ [°] | 2.12076778 ±5.8225900 | 2.12076778 ±5.8225900 | 2.12076778 ±5.82194309 | 2.12076778 ±5.8225900 | 2.12076778 ±5.8225900 |
| $\omega$ [°] | 34.686929715 ±4.098509955 | 34.686929715 ±4.098509955 | 34.686929715 ±4.098509955 | 34.686929715 ±4.098509955 | 34.686929715 ±4.098509955 |
| $r_4$ | **1.39482577632** **±0.052189395** | 0.954333337 ±0.010715192 | 0.95433333 ±0.010715192 | **1.39482577632** **±0.052189395** | 0.954333337 ±0.010715192 |
| $r_1$ | **0.0148487230** **±0.071517683** | 0.0101594275 ±0.04893072 | 0.0101594275 ±0.04893072 | **0.0148487230** **±0.071517683** | 0.0101594275 ±0.04893072 |
| $r_2$ | **-0.032969738** **±0.077595318** | -0.02255774 ±0.05308425 | -0.02255774 ±0.05308425 | **-0.032969738** **±0.077595318** | -0.02255774 ±0.05308425 |
| $r_3$ | **-0.4351354781** **±0.05222766** | -0.297717679 ±0.03411742 | -0.297717679 ±0.03411742 | **-0.4351354781** **±0.05222766** | -0.297717679 ±0.03411742 |
| $s_1$ | | | 75.0934 ±11.9697 | **51.3785932** **±9.36527** | 75.0934 ±11.9697 |
| $s_2$ | | | 80.9610 ±12.02106 | **55.3931683** **9.36228** | 80.9610 ±12.02106 |
| $s_3$ | | | -14.218104 ±19.7218 | **-9.727961** **±13.48190** | -14.218104 ±19.7218 |
| $s_4$ | | | -3.32126 ±7.23696 | **-2.272390** **±4.94820** | -3.32126 ±7.23696 |
| $\widehat{\sigma_0}[m]$ | 10.7709 | 10.7709 | 10.7709 | 10.7709 | 10.7709 |

Note: Different values in column (1) and column (4) are shown in bold. The differences in column (1) and column (4) are due to the use of scaled quaternions. The results of Zeng et al. (2024) column (3) and the present methods column (5) are based on unit quaternions.

It can be seen that the results of the QA-based 3D transformation we performed with unit quaternions column (2) are fully compatible with the DQA-based results with unit quaternions column (5). The computed residual of coordinates in the first and second systems by QA and DQA methods are shown in Table 9. The variance-covariance matrix of all parameters is given in Table 10-11. In summary, the transformation parameters and precision of the parameters given by the QA and DQA methods are exactly the same. However, depending on whether scaled or unit quaternions are used in the calculations, there will be different quaternion values and different quaternion precision. This is also quite natural.

**Table 9** Residuals

| Points | Residual[m] | | | | | |
|---|---|---|---|---|---|---|
| | $v_x$ | $v_y$ | $v_z$ | $v_X$ | $v_Y$ | $v_Z$ |
| 1 | 1.9534 | -1.6429 | -4.8511 | -0.4262 | 1.1391 | 2.2595 |
| 2 | 3.2523 | -7.7132 | 2.4255 | 0.8548 | 3.8425 | -1.0719 |
| 3 | -8.6615 | 1.8208 | -1.9404 | 2.8032 | -3.0124 | 1.0293 |
| 4 | 3.2989 | 3.1293 | 1.2128 | -2.0729 | -0.3233 | -0.6723 |

**Table 10** The part of $C$ variance-covariance matrix belonging to the ($\lambda$, $r_1$, $r_2$, $r_3$, $r_4$, $s_1$, $s_2$, $s_3$, $s_4$) parameters

| | | | | | | | | |
|---|---|---|---|---|---|---|---|---|
| 0.0233 | 0.0000 | 0.0000 | 0.0000 | 0.0000 | -1.0498 | -0.9073 | -0.1395 | -0.0538 |
| 0.0000 | 0.0024 | -0.0003 | 0.0000 | -0.0000 | 0.0096 | 0.0270 | -0.6107 | -0.3483 |
| -0.0000 | -0.0003 | 0.0028 | -0.0000 | 0.0001 | -0.0376 | -0.0007 | 0.8637 | 0.0582 |
| -0.0000 | 0.0000 | -0.0000 | 0.0012 | 0.0004 | 0.3023 | -0.3265 | -0.0146 | 0.0018 |
| 0.0000 | -0.0000 | 0.0001 | 0.0004 | 0.0001 | 0.0933 | -0.1021 | 0.0224 | 0.0056 |
| -1.0498 | 0.0096 | -0.0376 | 0.3023 | 0.0933 | 143.2756 | -43.8112 | -6.7913 | 2.5779 |
| -0.9073 | 0.0270 | -0.0007 | -0.3265 | -0.1021 | -43.8112 | 144.5059 | -0.0372 | -3.4099 |
| -0.1395 | -0.6107 | 0.8637 | -0.0146 | 0.0224 | -6.7913 | -0.0372 | 388.9484 | 96.0516 |
| -0.0538 | -0.3483 | 0.0582 | 0.0018 | 0.0056 | 2.5779 | -3.4099 | 96.0516 | 52.3736 |

**Table 11** The variance-covariance *Cpar* matrix ($\varepsilon, \psi, \omega, t_X, t_Y, t_Z$)

| | | | | | |
|---|---|---|---|---|---|
| 0.01054 | -0.001639 | -0.0002172 | -0.263 | -0.2182 | 2.495 |
| -0.001639 | 0.010330 | 0.0001593 | 0.3052 | -0.1286 | -1.735 |
| -0.000217 | 0.000159 | 0.005117 | -0.5107 | 1.1930 | -0.06824 |
| -0.263 | 0.3052 | -0.5107 | 410.9 | 0.8242 | -57.93 |
| -0.2182 | -0.1286 | 1.193 | 0.8242 | 405.2 | -12.61 |
| 2.495 | -1.735 | -0.06824 | -57.93 | -12.61 | 844.8 |

## Conclusion

All stages of DQA-based symmetric 3D coordinate transformation are explained in detail step by step. The constrained model used in the transformation is simpler and easier to understand than the TLS method. The second method we use is the simplified method (unconstrained model), which is much simpler. The two numerical applications, Uygur et al. (2022), Zeng et al. (2024) and presented two methods that gave the same transformation parameters and residual. The different values in the Zeng et al. (2024) method may be due to a calculation error or typing error. Numerical applications show that DQA-based 3D transformation results are fully compatible with QA-based results and that the opposite claim (that DQA-based results provide more accurate results) is invalid. When compared in terms of time, the result was surprisingly in favour of DQA. It was seen that the solution was reached in 7 iterations in DQA and 9 iterations in both of QA and 12 iterations in classical Euler Angles method.

It is noticeable that our two solution algorithms are much simpler and easier to understand compared to Zeng et al.'s (2024) method. It is evaluated that four contributions were made in the study; First, in the classical dual quaternion transformation algorithm, the number of equations to be written per control point is reduced from four to three. The total number of equations to be written was reduced from *4n+1* to *3n+2* and even to *3n* in the simplified model (*n* : number of control points). The simplified method reduces the number of unknowns from 9 to 7. Reducing the number of equations and unknowns will reduce the calculation load. By obtaining a solution by using the modified misclosure vector against the non-convergence problem in ill-conditioned and unstable equation systems, thirdly, by using the CEIV and EIV models, which are simpler and more understandable, against the TLS method. In addition, claims in the literature that DQA-based transformations produce more accurate results have been refuted.

Numerical applications have shown that, in addition to the constrained method, the simplified unconstrained solution model, which is characterised by its simplicity and ease, can be used without difficulty. It has also been tested that the new symmetric coordinate transformation based on the dual quaternion algorithm (constrained and unconstrained method) presented here can be successfully applied to all types of 3D transformation problems (large or small scale, large or small rotation angles). The fact that the two different models we proposed gave exactly the same results as Uygur et al. (2022) and Zeng et al. (2024) (except actual geodetic datum transformation case, Page 13, Table 6,7,8) demonstrated the validity of our methods. The proposed algorithm can perform both 2D and 3D symmetric similarity transformations. For the 2D transformation, it is sufficient to replace the $z$ and $Z$ coordinates in both systems with zero.

**Acknowledgements**

I would like to thank the editor and the anonymous reviewers and for providing valuable comments and suggestions that helped improve the manuscript greatly.

**Conflicts of Interest:** The authors declare no conflict of interest.

# References

Bektas S (2022) A new algorithm for 3D similarity transformation with dual quaternion, Arabian Journal of Geosciences,15:1273,https://doi.org/10.1007/s12517-022-10457-z

Bektas S (2024a) An expanded dual quaternion algorithm for 3D Helmert transformation and determination of the VCV matrix of the transformation's parameters, Journal of Spatial Science, https://doi: 10.1080/14498596.2023.2274997

Bektas S (2024b) Best (orthogonal) fitting ellipsoid with quaternions. Survey Review, 56(396), 249–264. https://doi.org/10.1080/00396265.2023.2225899

Bektas S (2024c) Comment on Ioannidou S; Pantazis G. Helmert Transformation Problem. From Euler Angles Method to Quaternion Algebra. ISPRS Int. J. Geo-Inf. 2020, 9, 494, ISPRS International Journal of Geo-Information 13, no. 10: 359. https://doi.org/10.3390/ijgi13100359

Clifford W.K (2007) Mathematical Papers. AMS Chelsea Publishing; American Mathematical Society: Providence, RI, USA.

Fang X (2015) Weighted total least-squares with constraints: a universal formula for geodetic symmetrical transformations. J Geod 89:459–469

Felus YA, Burtch RC (2009) On symmetrical three-dimensional datum conversion. GPS Solut 13:65–74

Grafarend E, and Awange J (2003) Nonlinear analysis of the three-dimensional datum transformation [conformal group C7(3)]. Journal of Geodesy, 77(1–2), 66–76. https://doi:10.1007/s00190-002- 0299-9

Hamilton W.R (1853) Lectures on Quaternions: Containing a Systematic Statement of a New Mathematical Method; University Press: Berlin, Germany.

Ioannidou S, Pantazis G (2020) Helmert Transformation Problem. From Euler Angles Method to Quaternion Algebra. ISPRS Int. J. Geo-Inf. 2020, 9, 494 14 of 14

Li B, Shen Y, Zhang X, Li C, Lou L (2013) Seamless multivariate affine error-invariables transformation and its application to map rectification. International Journal of Geographical Information Science 27(8):1572-1592

Mahboub V (2016). A weighted least-squares solution to a 3-D symmetrical similarity transformation without linearization. Stud Geophys Geod **60**, 195–209 . https://doi.org/10.1007/s11200-015-1109-1

Mercan H, Akyilmaz O, Aydin C (2018) Solution of the weighted symmetric similarity transformations based on quaternions. J Geod 92:1113–1130

Mikhail E.M, and Ackermann F (1976) Observations and least squares. New York: IEP–A Dun–Donnelley

Öztürk E, Şerbetçi M (1982) Dengeleme Hesabı Cilt III (Adjustment Calculus), Karadeniz Technical University Press, Yayın No:144, Trabzon, Turkey (in Turkish).

Uygur S Ö, Aydin C, Akyılmaz O (2022) Retrieval of Euler rotation angles from 3D similarity transformation based on quaternions, Journal of Spatial Science,67(2),255-272, https://doi: 10.1080/14498596.2020.1776170

Wang Y, Kun Y, Zheng N, Bian Z, Yang M (2023) A generalized weighted total least squares-based, iterative solution to the estimation of 3D similarity transformation parameters, Measurement,210 112563, https://doi.org/10.1016/j.measurement.2023.112563

Zeng H, Yi Q (2011) Quaternion-Based Iterative Solution of Three-Dimensional Coordinate Transformation Problem. J. Comput. 2011, 6, 1361–1368.

Zeng H, Fang X, Chang G, Yang R (2018) A dual quaternion algorithm of the Helmert transformation problem, Earth, Planets and Space 70:26 https://doi.org/10.1186/s40623-018-0792-x

Zeng H, Chang G, He H, Tu Y, Sun S, Wu, Y (2019) Iterative solution of Helmert transformation based on a unit dual quaternion, Acta Geodaetica et Geophysica 54:123–141,https://doi.org/10.1007/s40328-018-0241-0

Zeng H , Wang Z, Li J, Li S, Wang J, Li X (2024) Dual-quaternion-based iterative algorithm of the three dimensional coordinate transformation, Earth, Planets and Space (2024) 76:20 https://doi.org/10.1186/s40623-024-01967-z

Zhao Z, Li Z, Wang B (2024) A Novel Robust Quaternions-Based Algorithm for 3-D Symmetric Similarity Datum Transformation, IEEE Transactions on Instrumentation and Measurement, vol. 73, pp. 1-12, Art no. 1003012, https://doi: 10.1109/TIM.2024.3370773.